# AVERAGE OPTIMALITY FOR CONTINUOUS-TIME MARKOV DECISION PROCESSES IN POLISH SPACES[1]

By Xianping Guo and Ulrich Rieder

*Zhongshan University and Universitaet Ulm*

This paper is devoted to studying the average optimality in continuous-time Markov decision processes with fairly general state and action spaces. The criterion to be maximized is expected average rewards. The transition rates of underlying continuous-time jump Markov processes are allowed to be *unbounded*, and the reward rates may have *neither upper nor lower bounds*. We first provide *two* optimality inequalities with *opposed* directions, and also give suitable conditions under which the existence of solutions to the two optimality inequalities is ensured. Then, from the two optimality inequalities we prove the existence of optimal (deterministic) stationary policies by using the *Dynkin formula*. Moreover, we present a "semimartingale characterization" of an optimal stationary policy. Finally, we use a generalized Potlach process with control to illustrate the *difference* between our conditions and those in the previous literature, and then further apply our results to average optimal control problems of generalized birth–death systems, upwardly skip-free processes and two queueing systems. The approach developed in this paper is slightly different from the "optimality inequality approach" widely used in the previous literature.

**1. Introduction.** Continuous-time Markov decision processes (MDPs) have received considerable attention because many optimization models such as those in telecommunication and queueing systems are based on the processes involving continuous time. One of the most common optimality criterion in continuous-time MDPs is the *expected average criterion*, which has been studied by many authors. In this paper we are also concerned with this expected average criterion. As is well known, continuous-time MDPs can be

Received January 2004; revised November 2004.

[1]Supported by NSFC, NCET and the Foundation of Zhongshan University Advanced Research Center.

*AMS 2000 subject classifications.* 90C40, 93E20.

*Key words and phrases.* Average reward, general state space, optimality inequality, optimal stationary policy, semimartingale characterization.

specified by four *primitive data*: a state space $S$; an action space $A$ with subsets $A(x)$ of admissible actions, which may depend on the current state $x \in S$; transition rates $q(\cdot|x,a)$; and reward (or cost) rates $r(x,a)$. Using these terms, we now briefly describe some existing works on the expected average criterion. When the state space is *finite*, a bounded solution to the average optimality equation (AOE) and methods for computing optimal stationary policies have been investigated in [23, 26, 30]. Since then, most work has focused on the case of a *denumerable* state space; for instance, see [6, 24] for bounded transition and reward rates, [18, 27, 31, 34, 39, 41] for bounded transition rates but unbounded reward rates, [16, 35] for unbounded transition rates but bounded reward rates and [12, 13, 17] for unbounded transition and reward rates. For the case of an *arbitrary* state space, to the best of our knowledge, only Doshi [5] and Hernández-Lerma [19] have addressed this issue. They ensured the existence of optimal stationary policies. However, the treatments in [5] and [19] are restricted to *uniformly bounded* reward rates and *nonnegative* cost rates, respectively, and the AOE plays a key role in the proof of the existence of average optimal policies. Moreover, to establish the AOE, Doshi [5] needed the hypothesis that all admissible action sets are *finite* and the relative difference of the optimal discounted value function is *equicontinuous*, whereas in [19] the *assumption* about the existence of a solution to the AOE is imposed. On the other hand, it is worth mentioning that some of the conditions in [5, 19] are imposed on the *family* of weak infinitesimal operators deduced from all admissible policies, instead of the primitive data. In this paper we study the much more general case. That is, the reward rates may have *neither upper nor lower bounds*, all of the state and action spaces are *fairly general* and the transition rates are allowed to be *unbounded*. We first provide *two* optimality inequalities rather than *one* for the "optimality inequality approach" used in [16, 19], for instance. Under suitable assumptions we not only prove the existence of solutions to the two optimality inequalities, but also ensure the existence of optimal stationary policies by using the two inequalities and the *Dynkin formula*. Also, to verify our assumptions, we further give sufficient conditions which are imposed on the *primitive data*. Moreover, we present a semimartingale characterization of an optimal stationary policy. Finally, we use controlled generalized Potlach processes [4, 22] to show that all conditions in this paper are satisfied, whereas the earlier conditions fail to hold. Then we further apply our results to average optimal control problems of generalized birth–death systems and upwardly skip-free processes [1], a pair of controlled queues in tandem [28], and $M/M/N/0$ queue systems [25, 40]. It should be noted that, on the one hand, the optimality inequality approach used in the previous literature (see, e.g., [16, 19] for continuous-time MDPs and [20, 21, 31, 34] for discrete-time MDPs) is *not* applied to our case, because in our model the reward rates may have *neither upper nor lower bounds*. On the other hand, we not only replace the AOE with two optimality inequalities, but also relax the condition of the equicontinuity of the relative difference of optimal discounted value functions [5]. Therefore, the approach developed in this paper can be regarded as a modification of the optimality inequality approach widely used in the previous literature.

The rest of this paper is organized as follows. In Section 2 we introduce the optimal control problem. Our main results are given in Section 4 after some technical preliminaries in Section 3. We illustrated with examples our conditions and results in Section 5, and conclude in Section 6 with some general remarks.

## 2. The optimal control problem.

NOTATION. If $X$ is a Polish space (i.e., a complete and separable metric space), we denote by $\mathcal{B}(X)$ the Borel $\sigma$-algebra.

The model of continuous-time MDPs with which we are concerned is of the form

$$(2.1) \quad \{S, (A(x) \subset A), q(\cdot|x,a), r(x,a)\}$$

with the following components.

• The variable $S$ is the *state space*—a Polish space.

• The term $A(x)$ is a Borel subset of $A$ which denotes the set of admissible actions at state $x \in S$, where $A$ is the *action space*—a Polish space too. The set

$$(2.2) \quad K := \{(x,a) | x \in S, a \in A(x)\}$$

of pairs of states and actions is assumed to be a Borel subset of $S \times A$.

• The element $q(\cdot|x,a)$ in (2.1) denotes the *transition rates*, which satisfy the following properties for each $(x,a) \in K$ and $D \in \mathcal{B}(S)$:

$P_1$: The element $q(\cdot|x,a)$ is a *signed* measure on $\mathcal{B}(S)$, and $q(D|\cdot,\cdot)$ is Borel measurable on $K$.
$P_2$: For all $x \notin D \in \mathcal{B}(S)$, $0 \leq q(D|x,a) < \infty$,
$P_3$: There exists $q(S|x,a) = 0, 0 \leq -q(\{x\}|x,a) < \infty$.

Furthermore, the model is assumed to be *stable*, that is,

$$(2.3) \quad q(x) := \sup_{a \in A(x)} (-q(\{x\}|x,a)) < \infty \qquad \forall x \in S.$$

• The real-valued function $r(x,a)$ denotes the *reward rates* and it is assumed to be measurable on $K$. [Whereas $r(x,a)$ is allowed to take positive and negative values, it can be interpreted as a *cost rate* rather than a "reward" rate.]

We now define a randomized Markov policy.

DEFINITION 2.1 (Randomized Markov policies). Let $\Phi$ be the set of functions $\pi_t(B|x)$ on $[0,\infty) \times \mathcal{B}(A) \times S$ such that:

(1) For each $t \geq 0$, $\pi_t(\cdot|x)$ is a stochastic kernel on $A$ given $S$ such that $\pi_t(A(x)|x) = 1$ for all $x \in S$.

(2) For each $B \in \mathcal{B}(A)$ and $x \in S$, $\pi_t(B|x)$ is a Borel measurable function in $t \geq 0$.

A function $\pi_t(B|x)$ in $\Phi$ is called a *randomized Markov policy*. We will write $\pi_t(B|x)$ simply as $(\pi_t)$. The subscript "$t$" in $\pi_t$ indicates the possible dependence on time. A randomized Markov policy $\pi := (\pi_t) \in \Phi$ is called (deterministic) *stationary* if there exists a Borel measurable function $f$ on $S$ such that

$$f(x) \in A(x) \quad \text{and} \quad \pi_t(\{f(x)\}|x) = 1 \qquad \forall t \geq 0 \text{ and } x \in S.$$

For simplicity, we denote by $f$ this stationary policy $\pi$. The set of all stationary policies is denoted by $F$; this means that $F$ is the set of all measurable functions $f$ on $S$ with $f(x) \in A(x)$ for all $x \in S$. Obviously, $F \subset \Phi$.

By (1) above, w.l.o.g. we also regard $\pi_t(\cdot|x)$ as a probability measure on $A(x)$. Thus, for each fixed policy $\pi = (\pi_t) \in \Phi$, the associated transition rates $q(D|x, \pi_t)$ can be defined by

$$q(D|x, \pi_t) := \int_{A(x)} q(D|x, a) \pi_t(da|x)$$

(2.4)

for each $x \in S, D \in \mathcal{B}(S)$ and $t \geq 0$.

In particular, when $\pi$ is stationary (i.e., $\pi =: f \in F$), we write the left-hand side of (2.4) as $q(D|x, f(x))$. Then $q(D|x, \pi_t)$ is called an *infinitesimal generator* [for any fixed policy $\pi = (\pi_t) \in \Phi$]; see [5, 24], for instance. Its equivalent form can be found in [8]. As is well known, any (possibly substochastic and nonhomogeneous) *transition function* $\tilde{p}^\pi(s, x, t, D)$ that depends on $\pi$ such that

$$\lim_{\varepsilon \to 0^+} \frac{\tilde{p}^\pi(t, x, t+\varepsilon, D) - I_D(x)}{\varepsilon} = q(D|x, \pi_t)$$

for all $x \in S, D \in \mathcal{B}(S)$ and $t \geq 0$ is called a *Q-process* with transition rates $q(D|x, \pi_t)$, where $I_D(x)$ is the indication function of set $D$.

To guarantee the existence of such a $Q$-process, we need to introduce the class of admissible policies.

DEFINITION 2.2 (An admissible policy). A policy $(\pi_t)$ in $\Phi$ is said to be *admissible* if for each $x \in S$ the functions $\int_{A(x)} h(a) \pi_t(da|x)$ are continuous in $t \geq 0$ for all bounded measurable functions $h$ on $A(x)$. We denote by $\Pi$ the class of all admissible policies. Observe that $\Pi$ is nonempty because it contains $F$. Moreover, it is easy to provide an example for which $\Pi$ can be chosen to be *strictly larger* than $F$.

By $P_1$–$P_3$, (2.3), (2.4) and Definition 2.2, we have the following facts.

LEMMA 2.1. *For each $\pi := (\pi_t) \in \Pi$, the following statements hold.*

(a) *For each $x \in S, t \geq 0$ and $D \in \mathcal{B}(S)$:*

  ($a_1$) $q(D|x, \pi_t)$ *is a signed measure in $D \in \mathcal{B}(S)$;*
  ($a_2$) $0 \leq q(D|x, \pi_t) < \infty$ *when $x \notin D$;*
  ($a_3$) $q(S|x, \pi_t) = 0$, $0 \leq -q(\{x\}|x, \pi_t) < \infty$;
  ($a_4$) $q(D|x, \pi_t)$ *is continuous in $t \geq 0$ and measurable in $x \in S$.*

(b) *There exists a $Q$-process $\tilde{p}^\pi(s, x, t, D)$ with transition rates $q(D|x, \pi_t)$.*

PROOF. Parts ($a_1$)–($a_3$) follow from (2.4) and the definition of model (2.1), while part ($a_4$) follows from (2.3) and Definition 2.2. By (a) and Theorem 2 in [8], we see that (b) is also true. □

Lemma 2.1(b) guarantees the existence of a $Q$-process such as the *minimum $Q$-process* $p_{\min}^\pi(s, x, t, D)$ [i.e., $p_{\min}^\pi(s, x, t, D) \leq \tilde{p}^\pi(s, x, t, D)$ for any $Q$-process $\tilde{p}^\pi(s, x, t, D)$], which can be directly constructed from the transition rates $q(D|x, \pi_t)$; see [8, 11], for instance. However, as is known [4, 8], such a $Q$-process might not be regular; that is, we might have $p_{\min}^\pi(s, x, t, S) < 1$ for some $x \in S$ and $t \geq s \geq 0$.

To ensure the regularity of a $Q$-process, we use the "drift" conditions below.

ASSUMPTION A. There exist a measurable function $w \geq 1$ on $S$, constants $c > 0, b \geq 0$ and $M_q > 0$, such that:

(1) For all $(x, a) \in K$, and $\int_S w(y) q(dy|x, a) \leq -cw(x) + b$.
(2) For all $x \in S$, with $q(x)$ as in (2.3), $q(x) \leq M_q w(x)$.

REMARK 2.1. (a) For the case of uniformly bounded transition rates [i.e., $\sup_{x \in S} q(x) < \infty$], Assumption A(2) is not required because it is only used to guarantee the regularity of a $Q$-process.

(b) Assumption A(1) is used not only for the regularity of a (possibly nonhomogeneous) $Q$-process, but also for the finiteness of the expected average criterion (2.6) below. Moreover, Assumption A(1) is a variant of the "drift condition" (2.4) in [28] for *homogeneous $Q$-processes*.

Under Assumption A, by Theorem 3.2 in [11] we see that a $Q$-process with transition rates $q(D|x, \pi_t)$ is *regular*, that is, $p_{\min}^\pi(s, x, t, S) = 1$ for all $x \in S$ and $t \geq s \geq 0$. Thus, under Assumption A we write the regular $Q$-process $p_{\min}^\pi(s, x, t, D)$ simply as $p^\pi(s, x, t, D)$.

We now state the optimality problem with which we are concerned.

For a given (initial) distribution $\mu_s$ on $S$ at $s \geq 0$ and each fixed policy $\pi = (\pi_t) \in \Pi$, let $p^\pi(s, x, t, D)$ be the regular $Q$-process. Then, as in the proof in [15], we can show the existence of a unique probability space $(\Omega, \mathcal{B}(\Omega), \tilde{P}_\pi^{\mu_s})$ with $\Omega = (S \times A)^{[0,\infty)}$. Let $\xi(t)$ and $\eta(t)$ denote the state and action processes, respectively (i.e., the coordinate processes defined on $\Omega$), and let $\tilde{E}_\pi^{\mu_s}$ denote the expectation operator associated with $\tilde{P}_\pi^{\mu_s}$. We write $\tilde{P}_\pi^{s,x}$ for $\tilde{P}_\pi^{\mu_s}$ and $\tilde{E}_\pi^{s,x}$ for $\tilde{E}_\pi^{\mu_s}$ when $\mu_s$ is the Dirac measure at $x \in S$. Moreover, let

$$(2.5) \qquad r(x, \pi_t) := \int_{A(x)} r(x, a) \pi_t(da|x) \qquad \text{for all } x \in S \text{ and } t \geq 0.$$

We will write $r(x, \pi_t)$ as $r(x, f(x))$ when $\pi = (\pi_t) =: f \in F$. Then we have the following lemma.

LEMMA 2.2. *Suppose that Assumption A holds. Then, for each $x \in S$ and $\pi = (\pi_t) \in \Pi$:*

(a) *For all $t \geq s \geq 0$, $\tilde{E}_\pi^{s,x} r(\xi(t), \eta(t)) = \int_S r(y, \pi_t) p^\pi(s, x, t, dy)$.*
(b) *The element $\tilde{E}_\pi^{s,x} r(\xi(t), \eta(t))$ is Borel measurable in $t$ ($t \geq s \geq 0$).*

PROOF. Part (a) follows from a similar proof in [15], while part (b) follows from (a) and (2.5) because $p(s, x, t, D, \pi)$ is continuous in $t \geq s \geq 0$; see [8]. □

For each $x \in S$ and $\pi \in \Pi$, the *expected average criterion* $V(x, \pi)$ is defined as

$$(2.6) \qquad V(x, \pi) := \liminf_{T \to \infty} \frac{\int_0^T [\tilde{E}_\pi^{0,x} r(\xi(t), \eta(t))] \, dt}{T}.$$

DEFINITION 2.3. A policy $\pi^*$ in $\Pi$ is said to be (average) optimal if $V(x, \pi^*) \geq V(x, \pi)$ for all $\pi \in \Pi$ and $x \in S$.

The main goal of this paper is to give conditions for the existence of an optimal stationary policy.

For each $x \in S, s \geq 0$ and $\pi := (\pi_t) \in \Pi$, we denote by $E_{s,x}^\pi$ the expectation operator associated with the probability measure $P_{s,x}^\pi$ which is completely determined by $p^\pi(s, x, t, D)$. Then by pages 107–109 in [10] (or by Theorem 14.4, page 121 in [38] and the homogenization technique in [5]) there exists a Borel measurable Markov process $x(t)$ ($t \geq 0$) with values in $S$. Obviously, the so-called state process $x(t)$ is a continuous-time *jump* Markov process; its transition function is $p^\pi(s, x, t, D)$ determined by the transition rates $q(D|x, \pi_t)$.

Then by Lemma 2.2, (2.6) and (2.5) we have

$$(2.7) \qquad V(x, \pi) = \liminf_{T \to \infty} \frac{\int_0^T [E_{0,x}^\pi r(x(t), \pi_t)] \, dt}{T}.$$

Here, we understand that $x(t)$ is any sample path and these sample paths are distributed according to $P_{s,x}^\pi$. Hence, these sample paths have a dependence on $\pi$, $s$ and $x$. However, such dependence will be dropped for simplicity when there is no confusion.

By (2.7), w.l.o.g. we will limit ourselves to use $x(t)$ and the corresponding $P_{s,x}^\pi$ and $E_{s,x}^\pi$ throughout the following discussion. In particular, let $P_x^\pi := P_{0,x}^\pi$ and $E_x^\pi := E_{0,x}^\pi$.

**3. Preliminaries.** In this section, we give some preliminary lemmas that are needed to prove our main results.

LEMMA 3.1. *Suppose that Assumption A holds. Then, for each $\pi \in \Pi$,*

$$E_x^\pi w(x(t)) \leq e^{-ct} w(x) + \frac{b}{c} \qquad \forall \, x \in S \text{ and } t \geq 0$$

*with $w(x)$, $c$ and $b$ as in Assumption A.*

For the proof, see Theorem 3.1 in [11].

To prove our main results, in addition to the previous result, we also need some facts on the $\alpha$-discounted criterion defined by (3.1) below. For each discount factor $\alpha > 0$, $x \in S$ and $\pi \in \Pi$, the $\alpha$-discounted criterion $J_\alpha(x, \pi)$ and the corresponding optimal discounted value function $J_\alpha^*(x)$ are defined by

$$(3.1) \quad J_\alpha(x, \pi) := \int_0^\infty e^{-\alpha t} [E_x^\pi r(x(t), \pi_t)] \, dt \quad \text{and} \quad J_\alpha^*(x) := \sup_{\pi \in \Pi} J_\alpha(x, \pi),$$

respectively.

A policy $\pi^*$ in $\Pi$ is said to be $\alpha$-discounted optimal if $J_\alpha(x, \pi^*) = J_\alpha^*(x)$ for all $x \in S$.

To ensure the finiteness of both $V(x, \pi)$ and $J_\alpha(x, \pi)$, and the existence of $\alpha$-discounted optimal stationary policies, we give the following conditions.

ASSUMPTION B. (1) For each $x \in S$, $A(x)$ is compact.

(2) For each fixed $x \in S$, $r(x, a)$ is continuous in $a \in A(x)$, and the functions $\int_S u(y) q(dy|x, a)$ are continuous in $a \in A(x)$ for all bounded measurable functions $u$ on $S$ and also for $u := w$ as in Assumption A.

(3) For all $a \in A(x)$ and $x \in S$, $|r(x,a)| \leq Mw(x)$ with some constant $M > 0$.

(4) There exist a nonnegative measurable function $w'$ on $S$, and constants $c' > 0, b' \geq 0$ and $M' > 0$ such that [with $q(x)$ as in (2.3)]

$$q(x)w(x) \leq M'w'(x) \quad \text{and} \quad \int_S w'(y)q(dy|x,a) \leq c'w'(x) + b' \qquad \forall\, (x,a) \in K.$$

REMARK 3.1. Assumptions B(1) and B(2) are similar to the standard continuity–compactness hypotheses for discrete-time MDPs; see, for instance, [21, 31] and references therein. Under Assumptions A and B(3), by Lemma 3.1 we see that the values $V(x,\pi)$ and $J_\alpha(x,\pi)$ are both finite. Assumption B(4) allows us to use the *Dynkin formula*. On the other hand, if $q(x)$ or $r(x,a)$ is bounded, then Assumption B(4) is not required.

LEMMA 3.2. *Under Assumptions A and B, the following statements hold, with $\alpha > 0$.*

(a) *For all $x \in S$ and $\pi \in \Pi$, $|J_\alpha(x,\pi)| \leq \frac{M}{c+\alpha}w(x) + \frac{bM}{c\alpha}$.*

(b) *The optimal discounted value function $J_\alpha^*(x)$ satisfies the optimality equation*

$$(3.2) \qquad \alpha J_\alpha^*(x) = \sup_{a \in A(x)} \left\{ r(x,a) + \int_S J_\alpha^*(y) q(dy|x,a) \right\} \qquad \forall\, x \in S.$$

(c) *There exists an $\alpha$-discounted optimal stationary policy $f_\alpha^* \in F$.*

For the proof, see Theorem 3.3 in [11].

To state our final conditions, we need to introduce the concept of the weighted norm used in [20, 21]. For any fixed measurable function $h \geq 1$ on $S$, a function $u$ on $S$ is called $h$-bounded if the weighted norm of $u$, $\|u\|_h := \sup_{x \in S} \frac{|u(x)|}{h(x)}$, is finite. Such a function $h$ will be referred to as a *weight function*. We denote by $B_h(S)$ the Banach space of all $h$-bounded measurable functions $u$ on $S$.

ASSUMPTION C. There exist two functions $v_1, v_2 \in B_w(S)$ (with $w$ as in Assumption A) and some state $x_0 \in S$ such that

$$v_1(x) \leq h_\alpha(x) \leq v_2(x) \qquad \forall\, x \in S \text{ and } \alpha > 0,$$

where $h_\alpha(x) := J_\alpha^*(x) - J_\alpha^*(x_0)$ is the so-called relative *difference of the optimal discounted value function $J_\alpha^*(x)$*.

REMARK 3.2. (a) Assumption C is a variant of the conditions for discrete-time MDPs; see (SEN2) on page 132 in [34] and Assumption 5.4.1(b) in [20], for instance.

(b) It should be noted that the function $v_1$ in our Assumption C may *not* be bounded below, and so the $h_\alpha(x)$ may not be bounded below either. However, the corresponding $h_\alpha(x)$ in [20, 34] is assumed to be *bounded below*.

To verify Assumption C, we now provide some sufficient conditions.

LEMMA 3.3. *Under Assumptions A and B, each one of the following conditions* (a) *and* (b) *implies Assumption C.*

(a) *For each $f \in F$ there exists a probability measure $\mu_f$ on $\mathcal{B}(S)$ such that*

$$(3.3) \quad \left| E_x^f[u(x(t))] - \int_S u(y)\mu_f(dy) \right| \le Re^{-\rho t} w(x) \qquad \forall\, |u| \le w \text{ and } t \ge 0,$$

*where $R > 0$ and $\rho > 0$ are constants independent of $f$.*

(b) *For some integer $d \ge 1$, $S := [0, \infty)^d$ and $q(x)$ is locally bounded on $S$. Moreover, the following conditions are satisfied:*

($b_1$) *Drift condition. The function $w$ in Assumption A is nondecreasing in each component and, moreover,*

$$\int_S w(y) q(dy|x,a) \le -cw(x) + bI_{\{0^d\}}(x) \qquad \forall\, (x,a) \in K,$$

*where $0^d := (0, 0, \ldots, 0) \in S$.*

($b_2$) *Monotonicity condition. For each $x_k \in S, a_k \in A(x_k)$ ($k = 1, 2$) and monotone set $D$ [i.e., $I_D(x)$ is increasing in $x \in S$], if $x_1 \le x_2$ and $x_2 \notin D$, then $\tilde{q}(D|x_1, a_1) \le \tilde{q}(D|x_2, a_2)$, and $\tilde{q}(D^c|x_1, a_1) \ge \tilde{q}(D^c|x_2, a_2)$ when $x_1 \le x_2$ and $x_1 \in D$, where $\tilde{q}(D|x_k, a_k) := q(D|x_k, a_k) - q(\{x_k\}|x_k, a_k)I_D(x_k)$ and $D^c := S - D$ is the complement of set $D$.*

PROOF. (a) Whereas $w(x) \ge 1$ and $|r(x,a)| \le Mw(x)$ for all $x \in S$ and $a \in A(x)$, by Lemma 3.2(c) and (3.3) we have

$$|h_\alpha(x)| = \left| \int_0^\infty e^{-\alpha t}[E_x^{f^*_\alpha} r(x(t), f^*_\alpha(x(t))) - E_{x_0}^{f^*_\alpha} r(x(t), f^*_\alpha(x(t)))] \, dt \right|$$

$$\le MR \int_0^\infty e^{-(\alpha+\rho)t}[w(x) + w(x_0)] \, dt$$

$$\le \frac{MR}{\rho}[1 + w(x_0)]w(x) =: v_2(x),$$

which verifies Assumption C with $v_1(x) := -v_2(x)$.

(b) By Theorem 5.47 in [4], we see that for each $f \in F$ the corresponding Markov process $x(t)$ is stochastically ordered. Moreover, for each $x \in S, f \in F$ and $|u| \le w$, from the proof of (7.1) in [28] and condition (b), we have

$$\left| E_x^f[u(x(t))] - \int_S u(x)\mu_f(dx) \right| \le 2e^{-ct}\left[w(x) + \frac{c}{b}\right] \le 2\left(1 + \frac{c}{b}\right)e^{-ct}w(x),$$

which gives condition (a) and so Assumption C follows. □

Obviously, Lemma 3.3 is also true when $S = [0, \infty)^d$ and $I_{\{0^d\}}(x)$ in condition (b$_1$) are replaced with $S = [\beta_1, \infty) \times \cdots \times [\beta_d, \infty)$ and $I_{\{x^0\}}(x)$, respectively, where $x^0 := (\beta_1, \ldots, \beta_d) \in S, \beta_i \geq 0$ $(i = 1, \ldots, d)$.

The validity of conditions (a) and (b) in Lemma 3.3 can also be obtained in several ways. For instance, [28] uses Assumption A and monotonicity conditions. Other approaches that yield exponential ergodicity (3.3) can be seen in [3, 7, 36, 40], for instance.

To prove our main results by using the Dynkin formula, we need the following facts from Lemma 5.2 in [11].

LEMMA 3.4. *Suppose that Assumptions A and B hold. Take arbitrarily* $\pi := (\pi_t) \in \Pi$ *and* $x \in S$.

(a) *For each* $u \in B_{w+w'}$ *(with $w$ and $w'$ as in Assumptions A and B):*

(a$_1$) $\|E_x^\pi |u(x(t))|\|_{w+w'} \leq \frac{b+b'+c+c'}{c+c'} \|u\|_{w+w'} e^{(c+c')t}$ $\forall t \geq 0$;

(a$_2$) $\lim_{t \searrow s} e^{\int_s^v q(\{x\}|x, \pi_\delta) d\delta} E_{v,y}^\pi u(x(t)) = \lim_{t \searrow s} E_y^\pi u(x(t)) = u(y)$ *for all* $y \in S$ *and* $s \geq 0$.

(b) *For each* $u \in B_w(S)$ *and* $t \geq s \geq 0$:

(b$_1$) $L^\pi u(s, x) := \lim_{t \downarrow 0} t^{-1}[E_{s,x}^\pi u(x(s+t)) - u(x)] = \int_S u(y) q(dy|x, \pi_s)$;

(b$_2$) $\|E_{s,x}^\pi |L^\pi u(t, x(t))|\|_{w+w'} \leq \frac{\|u\|_w (c+c'+b+b'+2M')^2}{c+c'} e^{(c+c')(t-s)}$.

Lemma 3.4 shows that $L^\pi$ in Lemma 3.4(b) is the *extended generator* of the $Q$-process $p^\pi(s, x, t, D)$ and that the domain of $L^\pi$ contains $B_w(S)$.

Finally, for ease of reference, we state a "measurable selection theorem" from [20, 33].

LEMMA 3.5 (A measurable selection theorem). *Let $C(A)$ be the collection of all nonempty compact subsets of $A$, and let $D$ be a multifunction from $S$ to $C(A)$ such that $\bar{K} := \{(x, a) | x \in S, a \in D(x)\}$ is a Borel subset of $S \times A$. If $v(x, a)$ is a real-valued measurable function on $\bar{K}$ such that $v(x, a)$ is continuous in $a \in D(x)$ for each $x \in S$, then there exists a measurable function $f: S \to A$ such that $f(x) \in D(x)$ for all $x \in S$ and*

$$v(x, f(x)) = \max_{a \in D(x)} v(x, a) \quad \text{for each } x \in S.$$

*Moreover, the function $v^*(x) := \max_{a \in D(x)} v(x, a)$ is measurable in $x \in S$.*

Lemma 3.5 will be used to prove the existence of an optimal stationary policy.

**4. The main results and proof.** In this section, we prove our main results.

THEOREM 4.1. *Suppose that Assumptions* A, B *and* C *hold, and* $\pi = (\pi_t)$ *is in* $\Pi$.

(a) *If there exist a constant $g$ and a function $u \in B_w(S)$ such that*

$$g \geq r(x, \pi_t) + \int_S u(y) q(dy|x, \pi_t) \qquad \forall x \in S \text{ and } t \geq 0,$$

*then $g \geq V(x, \pi)$ for all $x \in S$.*

(b) *Similarly, if there exist a constant $g$ and a function $u \in B_w(S)$ such that*

$$g \leq r(x, \pi_t) + \int_S u(y) q(dy|x, \pi_t) \qquad \forall x \in S \text{ and } t \geq 0,$$

*then $g \leq V(x, \pi)$ for all $x \in S$.*

PROOF. (a) For each $x \in S$ and $T \geq 0$, under condition (a), by Lemma 3.4 and the Dynkin formula (page 141 or 146 in [9]), we have

$$
\begin{aligned}
E_x^\pi u(x(T)) - u(x) &= E_x^\pi \left[ \int_0^T L^\pi u(t, x(t))\, dt \right] \\
&\leq Tg - E_x^\pi \left[ \int_0^T r(x(t), \pi_t)\, dt \right] \\
&= Tg - \int_0^T E_x^\pi r(x(t), \pi_t)\, dt.
\end{aligned}
\tag{4.1}
$$

On the other hand, by Lemma 3.1 we have

$$E_x^\pi |u(x(T))| \leq \|u\|_w \left[ e^{-cT} w(x) + \frac{b}{c} \right],$$

which together with (4.1) and (2.7) gives (a).

(b) Similarly, we can prove (b). $\square$

THEOREM 4.2. *Suppose that Assumptions* A, B *and* C *hold.*

(a) *There exist a constant $g^*$, functions $u_1^*, u_2^* \in B_w(S)$ and a stationary policy $f^* \in F$ that satisfy the two optimality inequalities*

$$g^* \geq \sup_{a \in A(x)} \left\{ r(x,a) + \int_S u_1^*(y) q(dy|x,a) \right\} \qquad \forall x \in S; \tag{4.2}$$

$$g^* \leq \sup_{a \in A(x)} \left\{ r(x,a) + \int_S u_2^*(y) q(dy|x,a) \right\} \qquad \forall x \in S, \tag{4.3}$$

$$= r(x, f^*(x)) + \int_S u_2^*(y) q(dy|x, f^*(x)) \qquad \forall x \in S. \tag{4.4}$$

(b) *For all $x \in S$, $g^* = \sup_{\pi \in \Pi} V(x, \pi) = V(x, f^*)$.*

(c) *Any stationary policy $f \in F$ that realizes the maximum of* (4.3) *is optimal, and so $f^*$ in* (4.4) *is an optimal stationary policy.*

PROOF. (a) Let $x_0 \in S$ be as in Assumption C and let $\{\alpha_n\}$ be any sequence of discount factors such that $\alpha_n \to 0$ as $n \to \infty$. By Lemma 3.2(a), $|\alpha_n J^*_{\alpha_n}(x_0)|$ is bounded in $n \geq 1$. Therefore, there exist a subsequence $\{\alpha_k\}$ of $\{\alpha_n\}$ and a constant $g^*$ that satisfy

$$(4.5) \quad \lim_{k \to \infty} \alpha_k J^*_{\alpha_k}(x_0) = g^*, \qquad u_1^*(x) := \liminf_{k \to \infty} h_{\alpha_k}(x).$$

Since $|h_\alpha(x)| \leq |v_1(x)| + |v_2(x)|$ for all $x \in S$ and $\alpha > 0$ (by Assumption C), by (4.5) we have

$$u_1^* \in B_w(S) \quad \text{and} \quad \lim_{k \to \infty} \alpha_k h_{\alpha_k}(x) = 0 \qquad \forall x \in S.$$

On the other hand, take any real-valued measurable function $m$ on $S$ such that $m(x) > q(x) \geq 0$ for all $x \in S$. Then, for each $x \in S$ and $a \in A(x)$, by $P_1$–$P_3$, we see that $P(\cdot|x, a)$ defined by

$$(4.6) \quad P(D|x, a) := \frac{q(D|x, a)}{m(x)} + I_D(x) \qquad \text{for all } D \in \mathcal{B}(S)$$

is a *probability measure* on $\mathcal{B}(S)$. In fact, by $P_1$ we see that $P(D|x,a) := \frac{q(D|x,a)}{m(x)} + I_D(x)$ is completely additive in $D \in \mathcal{B}(S)$. Moreover, by $P_3$ we also see that $P(S|x,a) = 1$. Thus, it suffices to show that $0 \leq P(D|x,a) \leq 1$ for all $D \in \mathcal{B}(S)$. For $x \notin D$, by $P_1, P_2$ and $P_3$,

$$-q(\{x\}|x, a) = q(S - \{x\}|x, a)$$
$$= q(D|x, a) + q(S - D - \{x\}|x, a)$$
$$\geq q(D|x, a) \geq 0,$$

which together with (2.3) and $m(x) > q(x)$ gives that $P(D|x,a) = \frac{q(D|x,a)}{m(x)} \in [0, 1)$ (since $x \notin D$). When $x \in D$, it also follows from $P_1, P_2$ and $P_3$ that

$$-q(\{x\}|x, a) \geq -q(\{x\}|x, a) - q(D - \{x\}|x, a)$$
$$= -q(D|x, a) = q(S - D|x, a) \geq 0,$$

which together with (2.3) and $m(x) > q(x)$ yields that $-1 < \frac{q(D|x,a)}{m(x)} \leq 0$, and so $P(D|x,a) = \frac{q(D|x,a)}{m(x)} + 1 \in [0, 1)$ (since $x \in D$).

Noting that $h_\alpha(x) = J^*_\alpha(x) - J^*_\alpha(x_0)$, by (4.6) and $P_3$ we can rewrite (3.2) as

$$(4.7) \quad \frac{\alpha J^*_\alpha(x_0)}{m(x)} + \frac{\alpha h_\alpha(x)}{m(x)} + h_\alpha(x)$$
$$= \sup_{a \in A(x)} \left\{ \frac{r(x,a)}{m(x)} + \int_S h_\alpha(y) P(dy|x, a) \right\} \qquad \forall x \in S.$$

Thus, for each $k \geq 1$ and $x \in S$, by (4.7) we have

$$
\begin{aligned}
(4.8) \quad &\frac{\alpha_k J^*_{\alpha_k}(x_0)}{m(x)} + \frac{\alpha_k h_{\alpha_k}(x)}{m(x)} + h_{\alpha_k}(x) \\
&\geq \frac{r(x,a)}{m(x)} + \int_S h_{\alpha_k}(y) P(dy|x,a) \qquad \forall\, a \in A(x).
\end{aligned}
$$

Applying the extension of Fatou's lemma (8.3.7 in [21]), by (4.8) and (4.5) we get

$$\frac{g^*}{m(x)} + u_1^*(x) \geq \frac{r(x,a)}{m(x)} + \int_S u_1^*(y) P(dy|x,a) \qquad \forall\, x \in S \text{ and } a \in A(x).$$

This together with (4.6), yields

$$g^* \geq r(x,a) + \int_S u_1^*(y) q(dy|x,a) \qquad \forall\, x \in S \text{ and } a \in A(x),$$

which gives

$$(4.9) \qquad g^* \geq \sup_{a \in A(x)} \left\{ r(x,a) + \int_S u_1^*(y) q(dy|x,a) \right\} \qquad \forall\, x \in S,$$

and so (4.2) follows.

To prove (4.3), for each $x \in S$ and $k \geq 1$, let

$$(4.10) \quad u_2^*(x) := \limsup_{k \to \infty} h_{\alpha_k}(x), \qquad g_{\alpha_k}(x) := \sup\{h_{\alpha_m}(x) : m \geq k\}.$$

Then we have

$$(4.11) \quad u_2^* \in B_w(S), \qquad u_2^*(x) = \lim_{k \to \infty} g_{\alpha_k}(x) \quad \text{and} \quad g_{\alpha_k}(x) \geq h_{\alpha_k}(x),$$

which together with (4.7) gives

$$
\begin{aligned}
(4.12) \quad &\frac{\alpha J^*_{\alpha_k}(x_0)}{m(x)} + \frac{\alpha h_{\alpha_k}(x)}{m(x)} + h_{\alpha_k}(x) \\
&= \sup_{a \in A(x)} \left\{ \frac{r(x,a)}{m(x)} + \int_S h_{\alpha_k}(y) P(dy|x,a) \right\} \\
&\leq \sup_{a \in A(x)} \left\{ \frac{r(x,a)}{m(x)} + \int_S g_{\alpha_k}(y) P(dy|x,a) \right\}.
\end{aligned}
$$

Since $g_{\alpha_{k+1}} \leq g_{\alpha_k}$ for all $k \geq 1$, $\lim_{k \to \infty}[\sup_{a \in A(x)}\{\frac{r(x,a)}{m(x)} + \int_S g_{\alpha_k}(y) P(dy|x,a)\}]$ exists. Thus, by (4.5), (4.10) and (4.12) we have

$$
\begin{aligned}
(4.13) \quad &\frac{g^*}{m(x)} + u_2^*(x) \\
&\leq \lim_{k \to \infty}\left[\sup_{a \in A(x)} \left\{ \frac{r(x,a)}{m(x)} + \int_S g_{\alpha_k}(y) P(dy|x,a) \right\}\right] \qquad \forall\, x \in S.
\end{aligned}
$$

Also, for each fixed $x \in S$ and $k \geq 1$, by Assumption B, there exists $a_k(x) \in A(x)$ such that

$$
\begin{aligned}
(4.14) \quad &\sup_{a \in A(x)} \left\{ \frac{r(x,a)}{m(x)} + \int_S g_{\alpha_k}(y) P(dy|x,a) \right\} \\
&= \frac{r(x, a_k(x))}{m(x)} + \int_S g_{\alpha_k}(y) P(dy|x, a_k(x)).
\end{aligned}
$$

Since $A(x)$ is compact, there exists a subsequence $\{a_{k_i(x)}\}$ of $\{a_k(x)\}$ such that $\lim_{i \to \infty} a_{k_i}(x) =: a'(x) \in A(x)$. Noting that $\|g_{\alpha_k}\|_w \leq \|v_1\|_w + \|v_2\|_w$ for all $k \geq 1$, by (4.13), (4.14) and the extension of Fatou's lemma (8.3.7 in [21]) we obtain

$$
\begin{aligned}
\frac{g^*}{m(x)} + u_2^*(x) &\leq \lim_{i \to \infty} \left[ \sup_{a \in A(x)} \left\{ \frac{r(x,a)}{m(x)} + \int_S g_{\alpha_{k_i}}(y) P(dy|x,a) \right\} \right] \\
&= \lim_{i \to \infty} \left[ \frac{r(x, a_{k_i}(x))}{m(x)} + \int_S g_{\alpha_{k_i}}(y) P(dy|x, a_{k_i}(x)) \right] \\
&= \frac{r(x, a'(x))}{m(x)} + \int_S u_2^*(y) P(dy|x, a'(x)) \\
&\leq \sup_{a \in A(x)} \left\{ \frac{r(x,a)}{m(x)} + \int_S u_2^*(y) P(dy|x,a) \right\},
\end{aligned}
$$

which yields (4.3).

Moreover, by Assumption B and the extension of Fatou's lemma (8.3.7 in [21]), we see that $v(x,a) := r(x,a) + \int_S u_2^*(y) q(dy|x,a)$ is continuous in $a \in A(x)$ for each $x \in S$. Then, Lemma 3.5 with $D(x) := A(x)$ together with (4.3) gives the existence of $f^*$ $(\in F)$, satisfying (4.4). Thus, the proof of (a) is complete.

(b) For each $\pi = (\pi_t) \in \Pi$, from (4.2) we get

$$g^* \geq r(x,a) + \int_S u_1^*(y) q(dy|x,a) \qquad \forall\, a \in A(x) \text{ and } x \in S,$$

which together with (2.4) and (2.5), gives

$$g^* \geq r(x, \pi_t) + \int_S u_1^*(y) q(dy|x, \pi_t) \qquad \forall\, t \geq 0 \text{ and } x \in S.$$

Thus, by Theorem 4.1(a) with $u := u_1^*$, we have

$$g^* \geq V(x, \pi) \qquad \forall\, x \in S \text{ and } \pi \in \Pi,$$

and so

$$(4.15) \qquad g^* \geq \sup_{\pi \in \Pi} V(x, \pi) \qquad \forall\, x \in S.$$

Similarly, by (4.4) and Theorem 4.1(b) with $u = u_2^*$, we have

$$(4.16) \qquad g^* \leq V(x, f^*) \qquad \forall x \in S.$$

By (4.15) and (4.16) we have $g^* = V(x, f^*) = \sup_{\pi \in \Pi} V(x, \pi)$ for all $x \in S$, and so (b) follows.

(c) Obviously, (c) follows from the proof of (a) and (b). $\square$

REMARK 4.1. (a) From the proof of Theorem 4.2 we see that the approach used to prove Theorem 4.2 is *different* from the optimality inequality approach (e.g., [16, 19] for continuous-time MDPs and [20, 21, 31, 34] for discrete-time MDPs). In fact, there are *two* key steps in the proof of the existence of an (average) optimal stationary policy by using the optimality inequality approach. The first step is to obtain an inequality as in (4.15) by the *Abelian theorem* (e.g., [19, 37]), relating the average criterion $V(x, \pi)$ to the discounted criterion $J_\alpha(x, \pi)$. The other step is to get another inequality as in (4.16) from the *optimality inequality* as in (4.4). However, to use the Abelian theorem, the reward (or cost) rates have to be nonpositive (or nonnegative). Therefore, the optimality inequality approach in the previous literature is *not* applied to our case because in our model the reward rates may have *neither upper nor lower bounds*.

(b) From the proof of Theorem 4.2(a), we also see that properties $P_1$–$P_3$ about the transition rates play a particular role. In fact, without these properties, we can neither define the probability measure $P(\cdot|x, a)$ in (4.6) nor prove Theorem 4.2(a) by applying the extension of Fatou's lemma (8.3.7 in [21]) to the right-hand sides of (4.8) and (4.14). On the other hand, it does not seem to be possible to prove the existence of (average) optimal feedback policies in controlled stochastic differential equations (SDEs) (see [2, 19], for instance) by using the above approach because in controlled SDEs such properties $P_1$–$P_3$ fail to hold.

Theorem 4.2 ensures the existence of an (average) optimal stationary policy. We now gives an interesting semimartingale characterization of such a policy.

For each $x \in S, f \in F, u \in B_w(S)$ and any constant $g$, let

$$(4.17) \qquad \begin{aligned} \Delta(x; f, u, g) &:= r(x, f(x)) + \int_S u(y) q(dy|x, f(x)) - g, \\ \mathcal{F}_t &:= \sigma\{x(s) : 0 \leq s \leq t\} \end{aligned}$$

and define a continuous-time stochastic process

$$(4.18) \quad M_t(f, u, g) := \int_0^t r(x(s), f(x(s))) \, ds + u(x(t)) - tg \qquad \text{for each } t \geq 0.$$

THEOREM 4.3. *Suppose that Assumptions* A, B *and* C *hold.*

(a) *If $f^*$ is the optimal stationary policy obtained in Theorem 4.2, and $u_1^*, u_2^*$ and $g^*$ are from Theorem 4.2, then:*

(a$_1$) *For all $x \in S$, and $\{M_t(f^*, u_2^*, g^*), \mathcal{F}_t\}$ is a $P_x^{f^*}$-submartingale.*
(a$_2$) *For all $f \in F$ and $x \in S$, $\{M_t(f, u_1^*, g^*), \mathcal{F}_t\}$ is $P_x^f$-supermartingale.*

(b) *Conversely, if there exist a policy $\hat{f} \in F$, functions $u_1', u_2' \in B_w(S)$ and some constant $g$ such that:*

(b$_1$) $\{M_t(\hat{f}, u_2', g), \mathcal{F}_t\}$ *is a $P_x^{\hat{f}}$-submartingale for all $x \in S$ and*
(b$_2$) $\{M_t(f, u_1', g), \mathcal{F}_t\}$ *is $P_x^f$-supermartingale for all $f \in F$ and $x \in S$,*
*then the stationary policy $\hat{f}$ is (average) optimal.*

PROOF. For each $f \in F, u \in B_w(S), x \in S$ and constant $g$, we have

$$(4.19) \quad \begin{aligned} & E_x^f[M_t(f, u, g)|\mathcal{F}_s] \\ &= M_s(f, u, g) + E_x^f\left[\int_s^t \Delta(x(y); f, u, g)\, dy | \mathcal{F}_s\right] \end{aligned} \quad \forall t \geq s \geq 0.$$

In fact, from (4.17) and (4.18), we have

$$(4.20) \quad \begin{aligned} & E_x^f\left[\int_s^t \Delta(x(y); f, u, g)\, dy | \mathcal{F}_s\right] \\ &= E_x^f\left[\int_s^t r(x(y), f(x(y)))\, dy | \mathcal{F}_s\right] \\ &\quad + E_x^f\left[\int_s^t H(x(y); f, u)\, dy | \mathcal{F}_s\right] - (t-s)g, \end{aligned}$$

where $H(x; f, u) := \int_S u(y) q(dy|x, f(x))$. Using the Markov property, we obtain

$$E_x^f\left[\int_s^t H(x(y); f, u)\, dy | \mathcal{F}_s\right] = E_{x(s)}^f\left[\int_s^t H(x(y); f, u)\, dy\right],$$

which together with Lemma 3.4 and the Fubini's theorem gives

$$(4.21) \quad E_x^f\left[\int_s^t H(x(y); f, u)\, dy | \mathcal{F}_s\right] = \int_s^t [E_{x(s)}^f H(x(y); f, u)]\, dy.$$

Applying Lemma 3.4 and the Dynkin formula (e.g., page 141 or 146 in [9]), from (4.21) we obtain

$$(4.22) \quad E_x^f\left[\int_s^t H(x(y); f, u)\, dy | \mathcal{F}_s\right] = E_{x(s)}^f u(x(t)) - u(x(s)).$$

Thus, replacing (4.22) into (4.20) we get

$$(4.23) \quad \begin{aligned} & E_x^f\left[\int_s^t \Delta(x(y); f, u, g)\, dy | \mathcal{F}_s\right] \\ &= E_x^f\left[\int_s^t r(x(y), f(x(y)))\, dy | \mathcal{F}_s\right] \\ &\quad + E_{x(s)}^f u(x(t)) - u(x(s)) - (t-s)g. \end{aligned}$$

On the other hand, from (4.18) and the Markov property we have

$$E_x^f[M_t(f,u,g)|\mathcal{F}_s]$$
$$(4.24) \qquad = M_s(f,u,g) + E_x^f\left[\int_s^t r(x(y),f(x(y)))\,dy|\mathcal{F}_s\right]$$
$$- u(x(s)) + E_{x(s)}^f u(x(t)) - (t-s)g.$$

Finally, use (4.24) and (4.23) to obtain (4.19).

(a) For each $f \in F$ and $x \in S$, from (4.2) and (4.17) we have

$$\Delta(x; f, u_1^*, g^*) \leq 0,$$

which together with (4.19) implies that $\{M_t(f, u_1^*, g^*), \mathcal{F}_t\}$ is $P_x^f$-supermartingale. Similarly, we see that $\{M_t(f^*, u_2^*, g^*), \mathcal{F}_t\}$ is a $P_x^{f^*}$-submartingale and so (a) follows.

(b) For each $x \in S, f \in F$ and $u \in B_w(S)$, taking expectations in both sides of (4.19) gives

$$(4.25) \qquad E_x^f M_t(f,u,g) = E_x^f M_s(f,u,g) + E_x^f\left[\int_s^t \Delta(x(y); f, u, g)\,dt\right]$$
$$\forall t \geq s \geq 0.$$

By Lemma 3.4 and (4.17), $\Delta(\cdot; f, u, g)$ belongs to $B_{w+w'}(S)$. Thus, by condition (b$_1$) we have

$$E_x^{\hat{f}}[M_t(\hat{f}, u_2', g)] \geq E_x^{\hat{f}}[M_s(\hat{f}, u_2', g)] \qquad \forall t \geq s \geq 0.$$

Then, by (4.25) and Fubini's theorem, we get

$$E_x^{\hat{f}}\left[\int_s^t \Delta(x(y); \hat{f}, u_2', g)\,dy\right] = \int_s^t E_x^{\hat{f}}[\Delta(x(y); \hat{f}, u_2', g)]\,dy$$
$$\geq 0 \qquad \forall t \geq s \geq 0$$

and so

$$E_x^{\hat{f}} \Delta(x(t); \hat{f}, u_2', g) \geq 0 \qquad \forall \text{ a.e. } t \geq 0.$$

Therefore, there exists a sequence $t_n \downarrow 0$ as $n \to \infty$ such that

$$(4.26) \qquad E_x^{\hat{f}} \Delta(x(t_n); \hat{f}, u_2', g) \geq 0 \qquad \forall n \geq 0 \text{ and } x \in S.$$

Since $\Delta(\cdot; \hat{f}, u_2', g) \in B_{w+w'}(S)$, letting $n \to \infty$ in (4.26), from Lemma 3.4(a) we get

$$\Delta(x; \hat{f}, u_2', g) \geq 0$$

and so

$$(4.27) \qquad g \leq r(x, \hat{f}(x)) + \int_S u_2'(y) q(dy|x, \hat{f}(x)) \qquad \forall x \in S.$$

Then, by (4.27) and Theorem 4.1(b), we get

(4.28) $$g \leq V(x, \hat{f}) \qquad \forall\, x \in S.$$

Similarly, as in the proof of (4.27), by condition ($b_2$) we have

$$g \geq r(x, f(x)) + \int_S u_1'(y) q(dy|x, f(x)) \qquad \forall\, x \in S \text{ and } f \in F$$

and so

$$g \geq r(x, a) + \int_S u_1'(y) q(dy|x, a) \qquad \forall\, x \in S \text{ and } a \in A(x).$$

Then, by (2.3), (2.5) and Theorem 4.1(a) we have

(4.29) $$g \geq \sup_{\pi \in \Pi} V(x, \pi) \qquad \forall\, x \in S.$$

Combining (4.28) with (4.29) gives

$$V(x, \hat{f}) = \sup_{\pi \in \Pi} V(x, \pi)$$

and so (b) follows. □

Theorem 4.3 gives a semimartingale characterization of an optimal stationary policy.

**5. Examples.** In this section we will use five examples to illustrate our conditions and results.

EXAMPLE 5.1 (Optimal control of generalized Potlach processes in [4, 22]). The generalized Potlach process [4, 22] is a $Q$-process generated by the infinitesimal operator $L$ defined by (5.1) below. Here we are interested in the following optimal control problem.

Take $S := [1, \infty)^d$ with an integer $d \geq 1$. Then the generalized Potlach process can be generated by the operator $L$ defined by

(5.1) $$Lu(x, a_1) := \sum_{i=1}^d \int_0^\infty \left[ u\left( x - e_i x_i + y \sum_{j=1}^d p_{ij} x_i e_j \right) - u(x) \right] dF_\lambda(y)$$
for $x \in S$,

where $a_1 := (p_{ij})$ is a Markov transition matrix on $\{1, 2, \ldots, d\}$, $e_i$ is the $i$th unit vector in $\mathbf{R}^d$ and $F_\lambda(y)$ is a real-valued distribution function with a *parameter* $\lambda$, which can be regarded as a fixed reward fee. When the process is at state $x = (x_1, \ldots, x_d) \in S$, the cost incurred at each component $x_i$ is presented by $q_i \in [0, q_i^*]$, where $q_i^* > 0$ for all $i = 1, \ldots, d$. Let $a_2 := (q_1, \ldots, q_d)$. Here we interpret the parameters $a_1$ and $a_2$ as an *action* $a :=$ $(a_1, a_2)$, which belongs to a set $A_1 \times A_2$ of *available actions*. Suppose that $A_1$ is a *finite* set of Markov transition matrices $(p_{ij})$, $A_2 := [0, q_1^*] \times \cdots \times [0, q_d^*]$, $F_\lambda(y) := (1 - e^{-\lambda y})I_{[0,\infty)}(y)$ with $\lambda > 1$ and, for each $f \in F$, $x^{(1)}, x^{(2)} \in S$ such that $x^{(1)} \leq x^{(2)}$ with the semiorder, $x_i^{(1)} \sum_{j=1}^d p_{ij}^1 e_j \leq x_i^{(2)} \sum_{j=1}^d p_{ij}^2 e_j$ for all $(p_{ij}^1), (p_{ij}^2) \in A_1$ and $i = 1, \ldots, d$. For each $D \in \mathcal{B}(S), x \in S$ and $a_1 = (p_{ij}) \in A_1$, let

$$(5.2) \quad \tilde{q}(D|x, a_1) := \sum_{i=1}^d \int_0^\infty I_{D \setminus \{x\}}\left(x - e_i x_i + y \sum_{j=1}^d p_{ij} x_i e_j\right) \lambda e^{-\lambda y}\, dy.$$

Then, for each $x \in S$ and $a = (a_1, a_2) \in A := A_1 \times A_2$ with $a_1 := (p_{ij})$ and $a_2 := (q_1, \ldots, q_d)$, the transition rates $q(D|x, a)$ and the reward rates $r(x, a)$, which may depend on given parameter $\lambda$, are defined by

$$(5.3) \quad q(D|x, a) := \tilde{q}(D|x, a_1) - I_D(x)\tilde{q}(S|x, a_1)$$

and

$$(5.4) \quad r(x, a) := \sum_{i=1}^d \sum_{j=1}^d q_i p_{ij} x_j - \lambda(x_1 + \cdots + x_d),$$

respectively.

For each $x = (x_1, \ldots, x_d) \in S$, let $w(x) := x_1 + x_2 + \cdots + x_d$. Then, by (5.3) we have

$$(5.5) \quad q(x) := \sup_{a \in A(x)} [-q(\{x\}|x, a)] \leq d.$$

Moreover, by (5.2) and (5.3) we have

$$(5.6) \quad \int_S w(y) q(dy|x, a) \leq (x_1 + \cdots + x_d) \int_0^\infty \lambda(y - 1)e^{-\lambda y}\, dy = -\frac{(\lambda - 1)}{\lambda} w(x)$$

which together with (5.5) verifies Assumption A with $c := \frac{(\lambda-1)}{\lambda}$ and $b = 0$.

By (5.4) we have $|r(x, a)| \leq d(q_1^* + \cdots + q_d^* + \lambda)w(x)$ for all $x \in S$ and $a \in A$, which together with (5.5), Remark 3.1 and the finiteness of $A_1$, implies Assumption B.

Finally, we verify Assumption C. In fact, let $D$ be any monotone set in $S$ [i.e., $I_D(x)$ is increasing in $x$]. For each $f \in F$, $x^{(1)}, x^{(2)} \in S$ such that $x^{(1)} \leq x^{(2)}$. Let $f(x^{(1)}) =: ((p_{ij}^1), a_2^{(1)})$, $f(x^{(2)}) =: ((p_{ij}^2), a_2^{(2)})$. Then $x_i^{(1)} \sum_{j=1}^d p_{ij}^1 e_j \leq x_i^{(2)} \sum_{j=1}^d p_{ij}^2 e_j$. Thus, for each $i \in \{1, 2, \ldots, d\}$ and $y \geq 0$, we have

$$\xi^1(i, y, f(x^{(1)})) := (x^{(1)} - x_i^{(1)} e_i) + y x_i^{(1)} \sum_{j=1}^d p_{ij}^1 e_j$$

$$\leq (x^{(2)} - x_i^{(2)} e_i) + y x_i^{(2)} \sum_{j=1}^{d} p_{ij}^2 e_j$$

$$=: \xi^2(i, y, f(x^{(2)})),$$

which together with the monotonicity of set $D$, gives

$$\xi^2(i, y, f(x^{(2)})) \in D \quad \text{if } \xi^1(i, y, f(x^{(1)})) \in D$$

and

$$\xi^1(i, y, f(x^{(1)})) \notin D \quad \text{if } \xi^2(i, y, f(x^{(2)})) \notin D.$$

Thus, if $x^{(1)}, x^{(2)} \notin D$, by (5.2) we have

(5.7)
$$\begin{aligned}
\tilde{q}(D|x^{(1)}, f(x^{(1)})) &= \sum_{j=1}^{d} \int_0^\infty I_D(\xi^1(i, y, f(x^{(1)}))) \lambda e^{-\lambda y} \, dy \\
&\leq \sum_{j=1}^{d} \int_0^\infty I_D(\xi^2(i, y, f(x^{(2)}))) \lambda e^{-\lambda y} \, dy \\
&= \tilde{q}(D|x^{(2)}, f(x^{(2)})),
\end{aligned}$$

and if $x^{(1)}, x^{(2)} \in D$, we also have

(5.8)
$$\begin{aligned}
\tilde{q}(D^c|x^{(1)}, f(x^{(1)})) &= \sum_{j=1}^{d} \int_0^\infty I_{D^c}(\xi^1(i, y, f(x^{(1)}))) \lambda e^{-\lambda y} \, dy \\
&\geq \sum_{j=1}^{d} \int_0^\infty I_{D^c}(\xi^2(i, y, f(x^{(2)}))) \lambda e^{-\lambda y} \, dy \\
&= \tilde{q}(D^c|x^{(2)}, f(x^{(2)})),
\end{aligned}$$

and so it follows from Lemma 3.3(b) that Assumption C holds.

By the discussions above, we see that for Example 5.1 all conditions in this paper are satisfied. It should be noted that in Example 5.1 the state space is *not* denumerable and the reward rates have *neither upper nor lower bounds*; see (5.4). Therefore, the earlier conditions in [5, 6, 13, 16, 17, 18, 19, 23, 24, 26, 27, 30, 31, 34, 35, 39, 41] fail to hold because, except in [5, 19], the state spaces in the previous literature are all denumerable, while the reward rates in [5] and cost rates in [19] are uniformly bounded and bounded below, respectively.

EXAMPLE 5.2 (Optimal control of birth–death systems in [1, 4]). Consider a controlled birth–death system in which the state variable denotes a population size at any time $t \geq 0$. The birth rate is assumed to be a fixed constant $\lambda > 0$, but the death rates $\mu$ are assumed to be controlled by a decision-maker. Here we interpret any death rate $\mu$ as an *action* $a$ (i.e., $\mu =: a$). When the system's state is at $x \in S := \{0, 1, \ldots\}$, the decision-maker takes an action $a$ from a given set $A(x) \equiv [\mu_1, \mu_2]$ with $\mu_2 > \mu_1 > 0$, which increases or decreases the death rates given by (5.10) and (5.11) below. This action incurs a cost at rate $r_c(x, a)$. In addition, suppose that the benefit caused by each population is presented by $p > 0$ for each unit of time, and then the decision-maker gets a reward at rate $px$ for each unit of time during which the system remains in state $x$.

We now formulate this system as a continuous-time Markov decision process. The corresponding transition rates $q(y|x, a)$ are given as

(5.9) $\quad q(1|0, a) = -q(0|0, a) := \lambda \quad \forall a \in [\mu_1, \mu_2],$

(5.10) $\quad q(0|1, a) := a, \quad q(1|1, a) = -a - \lambda,$
$\quad q(2|1, a) := \lambda \quad \forall a \in [\mu_1, \mu_2].$

For each $x \geq 2$ and $a \in A(x) = [\mu_1, \mu_2]$,

(5.11) $\quad q(y|x, a) := \begin{cases} p_1 ax, & \text{if } y = x - 2, \\ p_2 ax, & \text{if } y = x - 1, \\ -(a + \lambda)x, & \text{if } y = x, \\ \lambda x, & \text{if } y = x + 1, \\ 0, & \text{otherwise,} \end{cases}$

where $p_1 \geq 0$ and $p_2 \geq 0$ are fixed constants and $p_1 + p_2 = 1$.

By the model's description we see that the reward rates $r(x, a)$ are of the form

(5.12) $\quad r(x, a) := px - r_c(x, a) \quad \text{for } (x, a) \in K := \{(x, a) : x \in S, a \in A(x)\}.$

We aim to find conditions that ensure the existence of an (average) optimal stationary policy. To do this, we consider the following assumptions:

$E_1$: There exists $\mu_1 - \lambda > 0$.
$E_2$: There exists $p_1 \leq \frac{\mu_1}{2\mu_2}$ with $p_1$ as in (5.11). (This condition obviously holds when $p_1 = 0$.)
$E_3$: The function $r_c(x, a)$ is continuous in $a \in A(x) = [\mu_1, \mu_2]$ for each fixed $x \in S$, and $c^*(x) := \sup_{a \in A(x)} |r_c(x, a)| < \tilde{M}(x + 1)$ for all $x \in S$ and some constant $\tilde{M} \geq 0$.

Under these conditions, we obtain the following.

PROPOSITION 5.1. *Under* $E_1$, $E_2$ *and* $E_3$, *the above controlled birth–death system satisfies Assumptions* A, B *and* C. *Therefore (by Theorem 4.2), there exists an optimal stationary policy.*

PROOF. We shall first verify Assumption A. Let $c := \frac{1}{2}(\mu_1 - \lambda) > 0$ (by $E_1$) and let $w(x) := x + 1$ for all $x \in S$. Then, from (5.9) and (5.10) we have

$$\text{(5.13)} \quad \sum_{y \in S} q(y|0, a) w(y) = \lambda \leq -cw(0) + \mu_1 + \lambda \quad \forall a \in A(x),$$

$$\text{(5.14)} \quad \sum_{y \in S} q(y|1, a) w(y) = -(a - \lambda) \leq -cw(1) \quad \forall a \in A(x).$$

Moreover, for each $x \geq 2$ and $a \in [\mu_1, \mu_2]$, from (5.11) we have

$$\text{(5.15)} \quad \begin{aligned} \sum_{y \in S} q(y|x, a) w(y) &= -(a + ap_1 - \lambda) x \\ &\leq -\tfrac{2}{3}(a + ap_1 - \lambda) w(x) \\ &\leq -cw(x). \end{aligned}$$

By (5.13)–(5.15) we have

$$\text{(5.16)} \quad \begin{aligned} \sum_{y \in S} q(y|x, a) w(y) &\leq -cw(x) + (\mu_1 + \lambda) I_{\{0\}}(x) \quad \forall a \in A(x) \text{ and } x \in S, \\ &\leq -cw(x) + \mu_1 + \lambda \quad \forall a \in A(x) \text{ and } x \in S, \end{aligned}$$

which gives Assumption A(1). On the other hand, by (5.9)–(5.11), we have $q(x) \leq (\mu_2 + \lambda)(x + 1) = (\mu_2 + \lambda) w(x)$ and so Assumption A(2) follows. Thus, Assumption A is true.

To verify Assumption B, by (5.12) and $E_3$ we have $|r(x, a)| \leq (p + \tilde{M}) w(x)$ for all $x \in S$. Thus, by (5.9)–(5.11) as well as $E_3$ we see that Assumptions B(1)–B(3) hold. To verify Assumption B(4), we let

$$w'(x) := (x + 1)(x + 2) \quad \text{for each } x \in S.$$

Then, by (5.9)–(5.11) we have

$$q(x) w(x) \leq (\mu_2 + \lambda) w'(x) \quad \forall x \in S,$$

$$\sum_{y \in S} q(y|x, a) w'(x) \leq 6\lambda w'(x) \quad \forall a \in [\mu_1, \mu_2] \text{ and } x \in S,$$

which imply Assumption B(4) with $M' := (\mu_2 + \lambda), c' := 6\lambda, b' := 0$.

Finally, we verify Assumption C. Since $0 \leq p_1 \leq \frac{\mu_1}{2\mu_2}$, by (5.9)–(5.11) we have that, for each fixed $f \in F$,

$$\sum_{y \geq k} q(y|x, f(x)) \leq \sum_{y \geq k} q(y|x + 1, f(x + 1)) \quad \forall x, k \in S \text{ such that } k \neq x + 1,$$

which together with Theorem 3.4 in [1], implies that the corresponding Markov process $x(t)$ is stochastically ordered. Thus, Assumption C follows from (5.16) and Lemma 3.3(b). □

EXAMPLE 5.3 (Optimal control of upwardly skip-free processes in [1]). The upwardly skip-free processes, also known as birth and death processes with *catastrophes*, belong to the category of *population processes* [1], Chapter 9, page 292, with the state space $S := \{0, 1, 2, \ldots\}$. Here we are interested in the average optimal control problem for such processes with catastrophes of *two* sizes, so the transition rates are of the form

$$(5.17) \quad q(y|x, a) := \begin{cases} \lambda x + a_1, & \text{if } y = x + 1, \\ -(\lambda x + \mu x + d(x, a_2) + a_1), & \text{if } y = x, \\ \mu x + d(x, a_2)\gamma_x^1, & \text{if } y = x - 1, \\ d(x, a_2)\gamma_x^2, & \text{if } y = x - 2, \\ 0, & \text{others,} \end{cases}$$

where $x \in S, a := (a_1, a_2)$, the constants $\lambda > 0, \mu > 0$, *immigration* rates $a_1 \geq 0$; $d(x, a_2)$ are nonnegative numbers that represent the rates at which the "catastrophes" occur and which are assumed to be controlled by decisions $a_2$ in some compact set $B(x)$, when the process is in state $x \geq 1$; the numbers $\gamma_x^1$ and $\gamma_x^2$ are nonnegative and such that $\gamma_x^1 + \gamma_x^2 = 1$ for all $x \geq 1$ and $\gamma_1^2 = 0$; and $\gamma_x^k$ is the probability that the process makes a transition to $x - k$ ($k = 1, 2$), given that a catastrophe occurs when the process is in state $x \geq 2$. For state $x = 0$, it is natural to let $d(0, a_2) \equiv 0$ and $\gamma_0^1 = \gamma_0^2 = 0$. On the other hand, we suppose that the immigration rates $a_1$ can also be controlled and so we interpret $a := (a_1, a_2)$ as an action. Thus, we may let the admissible action sets $A(0) := [0, b]$ and $A(x) := [0, b] \times B(x)$ for $x \geq 1$, with some constant $b > 0$. In addition, suppose that the damage caused by a catastrophe is represented by $p > 0$ for each unit of time and that it incurs a cost at rate $c(x, a_2)$ to take decision $a_2 \in B(x)$ at state $x \geq 1$. Let $c(0, \cdot) :\equiv 0$. Also, we assume that the benefits obtained by the transitions to $x - 1$ and $x - 2$ from $x$ ($\geq 2$) are represented by positive constants $q_1$ and $q_2$, respectively, and the benefit caused by each $a_1 \in [0, b]$ is represented by a real number $\tilde{r}(a_1)$. Then the reward rates are of the form

$$r(x, a) := \tilde{r}(a_1) - c(x, a_2) - pd(x, a_2) + q_1\gamma_x^1 d(x, a_2) + q_2\gamma_x^2 d(x, a_2)$$

for all $a = (a_1, a_2) \in A(x)$. As in the verifications of Assumptions A, B and C in Example 5.2, under the following conditions $F_1$–$F_3$, the above controlled upwardly skip-free processes satisfy Assumptions A, B and C, and, therefore (by Theorem 4.2), there exists an optimal stationary policy:

$F_1$: For all $x \geq 1$, $\mu - \lambda > 0$; $\gamma_{x+1}^2 \leq \inf_{\{a_2 \in B(x)\}} \frac{d(x, a_2) + \mu x}{d(x+1, a_2)}$.
$F_2$: There exists $b \leq \lambda - \mu + \inf_{\{x \geq 1, a_2 \in B(x)\}}\{d(x, a_2) + \gamma_x^2 d(x, a_2)\}$.
$F_3$: For each $x \in S$, the functions $\tilde{r}(a_1)$ and $c(x, a_2)$ are continuous in $(a_1, a_2) \in A(x)$, and $\sup_{a_2 \in B(x)} |d(x, a_2)| \leq L_1(x+1), \sup_{a_2 \in B(x)} |c(x, a_2)| < L_2(x+1)$ for some constants $L_1 > 0$ and $L_2 > 0$.

In particular, all of $F_1$–$F_3$ hold when $\lambda < \mu \leq b+\lambda$, $\tilde{r}(a_1) := \tau a_1$, $d(x, a_2) := 2a_2 x$, $\gamma_x^2 \leq \frac{1}{2} + \frac{\mu}{4\beta}$ and $B(x) := [b, \beta]$ for all $x \geq 1$, with some constants $\tau > 0$ and $\beta > b$.

EXAMPLE 5.4 (Optimal control of a pair of $M/M/1$ queues in tandem in [28]). Suppose that customers arrive as a Poisson stream with *unit* rate to the first queue, where they are serviced with mean service time $a_1^{-1}$. After service is completed at the first queue, each customer immediately departs and joins the second queue, where the mean service time is $a_2^{-1}$. After service is completed at the second queue, the customers leave the system with state space $S := \{0, 1, 2, \ldots\}^2$. Here, we interpret any given pair of mean service times $(a_1, a_2) =: a$ as an action and let corresponding action sets $A(x_1, x_2) \equiv [\mu_1, \mu_1^*] \times [\mu_2, \mu_2^*]$ with positive constants $\mu_1^* > \mu_1, \mu_2^* > \mu_2$. As in [28], let

$$w(x_1, x_2) := \sigma_1^{x_1-1} + \sigma_2^{x_1+x_2-1} + \gamma \sigma_1^{-\beta_1(x_1-1)} \sigma_2^{-\beta_2(x_1+x_2-1)},$$

where $\sigma_1 = 1.06, \sigma_2 = 1.03, \gamma = 0.4, \beta_1 = 1.5$ and $\beta_2 = 0.3$. Suppose that $\mu_1 \geq 3$ and $\mu_2 \geq 2$. Then, when $r(x_1, x_2, a)$ is *bounded* in all $(x_1, x_2, a)$ and *continuous* in $a \in A(x_1, x_2)$ for each fixed $(x_1, x_2) \in S$, from the argument in [28] and Lemma 3.3(b), we see that Assumptions A, B and C are all satisfied. In fact, under these parameter values, using the argument in [28] and a straightforward calculation, we can verify Assumption B and also Assumption A as well as the conditions (b) in Lemma 3.3 with $c := 0.002$ and $b := 0$. Therefore, there exists an average optimal stationary policy for this example with the above parameter values.

EXAMPLE 5.5 (Optimal control of $M/M/N/0$ queue systems in [25, 40]). Here the state space is $S := \{0, 1, 2, \ldots, N\}$ with some integer $N \geq 1$. Suppose that the arrival rate $\lambda$ is fixed but the service rates $\mu$ can be controlled. Therefore, we interpret service rates $\mu$ as actions, which may depend on the current states $x \in S$. We denote by $A(x)$ the action sets at state $x \in S$. Since there is no service in the queue at state 0, we may suppose that $A(0) := \{0\}$ for simplicity. Also, for each $x \geq 1$ we let $A(x) := [\mu_1, \mu_2]$ with constants $\mu_2 > \mu_1 > 0$. Then, the transition rates are given as $q(0|0, 0) = -\lambda = -q(1|0, 0)$ and $q(N|N, \mu) = -N\mu = -q(N-1|N, \mu)$ for all $\mu \in A(N)$. Moreover, for each $1 \leq x \leq N-1$ and $\mu \in A(x)$,

$$q(y|x, \mu) := \begin{cases} \lambda, & \text{if } y = x+1, \\ -(\lambda + \mu x), & \text{if } y = x, \\ \mu x, & \text{if } y = x-1, \\ 0, & \text{others.} \end{cases}$$

Thus, when $\mu_1 > \lambda$ and a reward rate function $r(x, \mu)$ is continuous in $\mu \in A(x)$ for all $x \in S$, as in the verification of Example 5.2, we see that this controlled $M/M/N/0$ queue system satisfies Assumptions A, B and C.

CONTINUOUS-TIME MARKOV DECISION PROCESSES 25

REMARK 5.1. In the verifications of Assumptions A, B and C for the five examples, a key step is the verification of Assumption C by using Lemma 3.3(b). This is due to the advantage that the drift and monotonicity conditions of Lemma 3.3(b) are imposed on the primitive data of the model. Here, we should note that these conditions have to be *uniform* with respect to the actions. In fact, such uniformity is used to show that the exponential convergence rate $\rho$ and the constant $R$ in (3.3) are *independent* of all stationary policies. On the other hand, other examples and approaches that yield exponential ergodicity (3.3) can be seen in [3, 7, 28, 36, 40], for instance. Finally, be warned that all of the underlying processes in this paper are continuous-time *jump* Markov processes, which can be determined by given transition rates (2.4) with the properties $P_1$–$P_3$.

**6. Concluding remarks.** In the previous sections we have studied the average optimality problem for continuous-time Markov decision processes in Polish spaces. Under suitable assumptions we have shown the existence of an optimal stationary policy. The approach developed to prove the existence of optimal stationary policies is different from the optimality inequality approach widely used in the previous literature. In addition, we have presented a semimartingale characterization for an optimal stationary policy. On the other hand, we believe that our formulation and approach are sufficiently general and, thus, provide a way to analyze other important problems, such as the problems of *bias optimality*, *Blackwell optimality* and *stochastic games* with average payoffs, which as far as we can tell have not been previously studied for continuous-time jump Markov processes with Polish spaces and unbounded transition rates. Research on these topics is in progress.

To conclude, it is worth noting that under our present conditions we cannot establish the average optimality equation by using the usual *diagonal argument*, because the state space may *not* be denumerable. We will give *additional* conditions under which the average optimality equation also holds in an upcoming paper.

**Acknowledgments.** We are very grateful to an Associate Editor and the anonymous referees for many fine comments and suggestions that have improved this paper.## REFERENCES

[1] ANDERSON, W. J. (1991). *Continuous-Time Markov Chains*. Springer, New York. MR1118840
[2] BOKAR, V. (1989). *Optimal Control of Diffusion Processes*. Longman Sci. Tech., Harlow. MR1005532
[3] CHEN, M. F. (2000). Equivalence of exponential ergodicity and $L^2$-exponential convergence for Markov chains. *Stochastic Process. Appl.* **87** 281–279. MR1757116

SCHOOL OF MATHEMATICS
AND COMPUTATIONAL SCIENCE
ZHONGSHAN UNIVERSITY
GUANGZHOU
PEOPLES REPUBLIC OF CHINA
E-MAIL: mcsgxp@zsu.edu.cn

ABTEILUNG MATHEMATIK VII
UNIVERSITAET ULM
ULM
GERMANY
E-MAIL: rieder@mathematik.uni-ulm.de